\documentclass[11pt]{article}

\usepackage[margin=1.1in]{geometry}

\usepackage{amsmath,amssymb,amsfonts}
\usepackage{algorithmic}
\usepackage{graphicx}
\usepackage{textcomp}
\usepackage{xcolor}
\def\BibTeX{{\rm B\kern-.05em{\sc i\kern-.025em b}\kern-.08em
    T\kern-.1667em\lower.7ex\hbox{E}\kern-.125emX}}

\newcommand{\bff}{{\bf f}}
\newcommand{\bfb}{{\bf b}}
\newcommand{\bfc}{{\bf c}}

\newcommand{\x}{{\bf x}}
\newcommand{\y}{{\bf y}}

\newcommand{\zero}{{\bf 0}}
\newcommand{\cc}{{\mathbb C}}

\newtheorem{theorem}{Theorem}[section]

\begin{document}

\title{GPU Accelerated Newton for Taylor Series Solutions \\
       of Polynomial Homotopies in Multiple Double Precision\thanks{Supported
by the National Science Foundation under grant DMS 1854513.}}

\author{Jan Verschelde\thanks{University of Illinois at Chicago,
Department of Mathematics, Statistics, and Computer Science,
851 S. Morgan St. (m/c 249), Chicago, IL 60607-7045
Email: {\tt janv@uic.edu}, URL: {\tt http://www.math.uic.edu/$\sim$jan}.}}

\date{18 June 2024}

\maketitle

\maketitle              
\begin{abstract}
A polynomial homotopy is a family of polynomial systems,
typically in one parameter $t$.
Our problem is to compute power series expansions of the 
coordinates of the solutions in the parameter $t$,
accurately, using multiple double arithmetic.
One application of this problem is the location of the nearest singular
solution in a polynomial homotopy, via the theorem of Fabry.
Power series serve as input to construct Pad\'{e} approximations.

Exploiting the massive parallelism of Graphics Processing Units (GPUs)
capable of performing several trillions floating-point operations per second,
the objective is to compensate for the cost overhead caused by arithmetic 
with power series in multiple double precision. 
The application of Newton's method for this problem requires the evaluation
and differentiation of polynomials, followed by solving a blocked lower
triangular linear system.
Experimental results are obtained on NVIDIA GPUs,
in particular the RTX 2080, RTX 4080, P100, V100, and A100.

Code generated by the CAMPARY software is used to obtain results in double
double, quad double, and octo double precision.  The programs in this study
are self contained, available in a public github repository 
under the GPL-v3.0 License.

\noindent {\bf Keywords and phrases.}
Graphics Processing Unit (GPU) \and multiple double arithmetic \and
Newton's method \and numerical analytic continuation \and Taylor series.
\end{abstract}
\section{Introduction}

Many problems in science and engineering require the solving
of a system of polynomial equations in several variables.
Writing polynomials backwards, starting at the lowest degree terms,
corresponds to power series truncated at the degree of the polynomials.
Homotopy methods~\cite{AY78}
define families of polynomial systems which connect 
a system that must be solved to a system with known solutions.
Continuation methods~\cite{AG03,Mor09}
track the solution paths defined by the homotopy.

This paper is another next step to accelerate
a new robust path tracking algorithm~\cite{TVV20a},
applying the linearization of~\cite{BV18}
and extending the multicore implementation of~\cite{TVV20b}.
As the robust path tracker of~\cite{TVV20a} scales well
to track millions of paths without error,
the context of this research effort is to scale the number of equations
and variables of the systems.
All code used in the experiments is publicly available
in the github repository of PHCpack~\cite{Ver99},
released under the GNU GPL license.

With multiple double arithmetic Taylor series developments for the
solution curves defined by polynomial homotopies can be computed accurately
using Newton's method.
The need for multiprecision in analytic continuation can be traced
back to~\cite{Hen66} via the direct quote
``Some reflection shows that in order to get a convergent process
{\em the early vectors $A_n^{(k)}$ (early with respect to $k$)
must be computed more accurately than the late ones}''
where the italics appear as in~\cite{Hen66}.
The early vectors refer to the earlier coefficients of the series.
This quote appeared in abbreviated form in the recent paper~\cite{Tre20}.

The Taylor series coefficients are input to algorithms
to construct Pad\'{e} approximants~\cite{BG96},
which are related to extrapolation methods~\cite{BRZ91}, \cite{Sid03}
and approximation algorithms~\cite{Tre13}.
An example of an application to electrical engineering is the
holomorphic embedding load flow method to solve power flow problems.
In the convergence study of~\cite{DLLLTSW21},
results are computed with 400 digits of precision, using~\cite{Advanpix}.

In this paper, the acceleration with Graphics Processing Units (GPUs) 
is applied to compensate for the computational overhead caused 
by the multiple double arithmetic 
needed to compute Taylor series accurately.
For this computation, Newton's method is accelerated using
kernels to differentiate polynomials at power series~\cite{Ver21}
and kernels for linear algebra operations~\cite{Ver22}.
Those kernels perform very well when given inputs of the proper size.
In a combined application, the dimension of the linear algebra
problem matches the number of variables of the polynomials,
which has an impact on the overall performance.

\subsection{Problem Statement}

The two main concerns are performance and convergence.
For performance, the input must be sufficiently large,
but still well conditioned enough to allow for Newton's method to converge.
In examining the scalability we address the first key question:
how much of the overhead can be compensated by GPU acceleration?
Prior work showed that teraflop performance was achieved in the
acceleration of the convolutions to evaluate and differentiate 
polynomials in several variables at power series~\cite{Ver21}
and in the acceleration of the blocked Householder QR~\cite{Ver22}
to solve linear systems in the least squares sense
in multiple double precision.
The differentiation adapts the reverse mode of algorithmic
differentiation~\cite{GW08} with power series arithmetic.
The acceleration of the blocked Householder QR~\cite{BV87}
is explained in~\cite{KCR09}, 
and also addressed in~\cite{BDT08} and~\cite{VD08}.
Additional related work on multicore and accelerated QR
can be found in~\cite{ABDK11}, \cite{AADF11},
\cite{KNDD12}, \cite{TDB10}, and~\cite{TNLD10}.
The acceleration of the back substitution algorithm
applies the formulas of~\cite{Hel78}, developed further
using the ideas of~\cite{NM01}.

The second question concerns the combination of various kernels
in the linearization of the power series.
Of the various different types of kernels that are launched,
which types require the most amount of time?

\subsection{Multiprecision Arithmetic}

A multiple double number is an unevaluated sum 
of nonoverlapping doubles.  
The renormalization and arithmetical operations are explained
in~\cite{MBDJJLMRT18}.

MPLAPACK~\cite{Mah22} supports quad double arithmetic and provides 
implementations of arbitrary multiprecision linear algebra operations.
A recent application to matrix-matrix multiplication is in~\cite{UK23}.
In~\cite{Kou23}, the 3M and 4M methods for complex multiplication
are applied in linear algebra operations on complex matrices.

In this paper,
for multiple double precision, the software libraries QDlib~\cite{HLB01}
and CAMPARY~\cite{JMPT16} are applied, customized as follows.
Instead of working with an array of double double numbers,
two arrays of doubles are used:
the first array for the most significant doubles
and the second one for the least significant doubles,
as this memory layout benefits memory coalescing,
especially for complex quad double and octo double numbers.
The GPU version of the QDlib library~\cite{LHL10}
uses the {\tt double2} and {\tt double4} types of the CUDA SDK,
which promote good memory access for double doubles and quad doubles,
but are not longer suited for complex quad doubles or octo doubles.

Multiple double precision is not true multiprecision in the sense
that one cannot select any number of bits for the fraction.
The other drawback is the limited size of the exponents
(limited to the 11 bits of the 64-bit hardware doubles),
which will prohibit the computation with infinitesimal values.
In the context of GPU acceleration, recent work of~\cite{IK20}
makes an interesting comparison with double double arithmetic:
``The double double arithmetic of CAMPARY performs best for the problem 
of matrix-vector multiplication.'' 
Concerning quad double precision, the authors of~\cite{IK20} write:
{\em ``the CAMPARY library is faster than our implementation;
however as the precision increases the execution time of
CAMPARY also increases significantly.''}
The advantage of multiple double arithmetic is that simple counts of
the number of floating-point operations quantify the cost overhead
precisely and the flops metrics for performance are directly applicable.

The goal of applying GPU acceleration to offset the cost of
multiple double precision arithmetic is related to the recent trend of
mixed precision~\cite{AABetal21},~\cite{HM22},~\cite{Kel22}.

Multiprecision is often deemed meaningless on problems with noisy
input data, but the problems of interest are solution paths 
highly sensitive to small changes in the continuation parameter, 
corresponding to Taylor series with a small convergence radius.
Reconditioning of the Taylor series is proposed in~\cite{VV22}.

\subsection{Numerical Condition of Taylor Series}

The convergence concern is closely related to
the numerical conditioning of the problem of power series solutions.
Consider the following classical result,
applied in~\cite{TVV20a} to detect nearby singularities.

\begin{theorem} 
[\mbox{\rm the ratio theorem of Fabry~\cite{Fab1896}}] \label{thmFabry}

\noindent {\em If for the series
$x(t) = c_0 + c_1 t + c_2 t^2 + \cdots + c_d t^d + c_{d+1} t^{d+1} + \cdots$,
we have
   $\displaystyle \lim_{d \rightarrow \infty} c_d/c_{d+1} = z$,
then
\begin{itemize}
\item $z$ is a singular point of the series, and 
\item it lies on the boundary of the circle of convergence of the series.
\end{itemize}
Then the radius of the disk of convergence is~$|z|$.}
\end{theorem}

In a numerical interpretation of this theorem,
observe that the smaller convergence radius,
the larger the growth in the coefficients,
because then $|c_{d+1}| > |c_d|$.
For example if $|z|$ is 1/2,
then for sufficiently large $d$, $|c_{d+1}| \approx 2 |c_d|$.
Thus, for series of order 64, we could observe
coefficients of magnitude $2^{64} \approx 1.8 \times 10^{19}$.
Therefore, for convergence, it is best to generate examples
which have a unique power series solution and use decaying coefficients
as in the series developments of exponential functions.

In this paper Newton's method is applied to compute Taylor series.
An alternative is to apply Fourier series methods as done in
analytic continuation, see e.g.~\cite{For81}.
These methods are very sensitive to a good choice of the step size 
for taking samples of the function to be differentiated.
As explained in~\cite{For81} a smaller step size benefits the
lower order coefficients while a larger step size may be needed
to compute the higher order coefficients accurately.

The main experimental result of this research is that
on systems of 1024 equations in 1024 variables,
when doubling the precision from quad double to octo double,
the increase in wall clock time is much less than what can
be predicted from the cost overhead factors of the multiple
double arithmetical operations as the increase in the time
occupied by the kernels is significant.
As expected, with multiple double arithmetic, 
the arithmetic intensity of the computations increases
and the GPU becomes fully occupied.

The next section contains a high level description of Newton's method
on power series, using linearizations.  The setup of the test problems
with the levels of precision is justified in the third section.
As defined in the fourth section,
Newton's method is executed using a staggered progression of the
order of the power series.  To predict the performance of
the accelerated code, in section five the arithmetic intensities
of different kernels are computed.  Section six contains the
results of the computational experiments.

\section{Linearized Series and Newton's Method}

Instead of working with vectors and matrices of power series,
we work with series that have as coefficients vectors and matrices.
For example, for series of order 4, we solve
$A(t) \x(t) = \bfb(t)$, with $n$-by-$n$ matrices $A_0$, $A_1$, $A_2$, $A_3$:
\begin{eqnarray}
    A(t) & = & A_0 + A_1 t + A_2 t^2 + A_3 t^3, \\
   \x(t) & = & \x_0 + \x_1 t + \x_2 t^2 + \x_3 t^3, \\
   \bfb(t) & = &\bfb_0 + \bfb_1 t + \bfb_2 t^2 + \bfb_3 t^3.
\end{eqnarray}
The linearization of $A(t) \x(t) = \bfb(t)$ leads to
\begin{equation} \label{eqToeplitz}
  \left[
     \begin{array}{cccc}
        A_0 & & &  \\
        A_1 & A_0 & &  \\
        A_2 & A_1 & A_0 &  \\
        A_3 & A_2 & A_1 & A_0  \\
     \end{array}
  \right]
  \left[
     \begin{array}{c}
        \x_0 \\ \x_1 \\ \x_2 \\ \x_3
     \end{array}
  \right]
  =
  \left[
     \begin{array}{c}
        \bfb_0 \\ \bfb_1 \\ \bfb_2 \\ \bfb_3
     \end{array}
  \right].
\end{equation}
Even as least squares solutions provide an accuracy close to machine
precision, in the forward substitution of the solution of the
lower triangular block Toeplitz system in~(\ref{eqToeplitz})
the errors in the coefficient vectors $\x_k$ of $\x(t)$ propagate
through the updates of the right hand side vectors~$\bfb_k$.
If we lose two decimal places of accuracy in each step,
then for series of order four up to eight decimal places
may be lost in the last coefficient vector~$\x_3$.
A more extensive error analysis was made in~\cite{TVV20b}.

Newton's method takes as input a system of polynomials in several
variables, with power series truncated to the same degree
and produces a sequence of power series.
As an operator, this version of Newton's method can be considered
as turning a problem in many variables where all variables are
interdependent into a sequence of power series for each separate
variable, thus removing the interdependencies among the variables.
A high level description of Newton's method is shown 
in the pseudo code below.

\begin{center}
\begin{tabbing}
$\quad$ \= Input: \= $\bff(\x(t))$, 
          system with power series coefficients; \\
\>      \> $\x_0 = \x(t)$,
          initial leading coefficients; \\
\>      \> $N$, the maximum number of iterations; \\
\>      \> $\epsilon$, the tolerance on the accuracy. \\
\> Output: $i$, the number of iterations; \\
\> \hphantom{Output:}
        $\x(t)$, if $i \leq N$, then $\| \bff(\x(t)) \| < \epsilon$. \\
\> $\x(0) := \x_0$ \\
\> for~\= $i$ from 1 to $N$ do \\
\>     \> $\displaystyle A(t), \x(t) := \partial ( \bff(\x(t)) )$ \\
\>     \> exit when $\| \bfb(t) \| < \epsilon$ \\
\>     \> $\Delta \x(t) := A(t) \backslash \bfb(t)$ \\
\>     \> report $\| \bfb(t) - A(t) \Delta \x(t) \|$ \\
\>     \> $\x(t) := \x(t) + \Delta \x(t)$ \\
\>     \> $\bfb(t) := \bfb(t) -  A(t) \x(t)$ \\
\>     \> exit when $\| \Delta \x(t) \| < \epsilon$ 
\end{tabbing}
\end{center}

The two computationally intensive operations
are the evaluation and differentiation $\partial$ and the linear system
solving $\backslash$.  
The residual computation $\| \bfb(t) - A(t) \Delta \x(t) \|$
is for monitoring the convergence of the method and can be omitted.

The high level description does not incorporate
the staggered nature of the updates.
As will be explained in section~\ref{secStaggered} below
not all coefficient vectors of the series are involved
in all stages.

\section{Columns of Monomials}

We start by considering polynomial systems
of the form $\x^E = \bfb(t)$, where $E$ is an exponent matrix
which contains in its rows the exponent vectors of the variables $\x$.
The right hand side $\bfb(t)$ is a vector of power series.
For example, $n=3$, $\x = [x_1, x_2, x_3]$:
\begin{equation}
  E =
  \left[
     \begin{array}{ccc}
        1 & 0 & 0 \\
        1 & 1 & 0 \\
        1 & 1 & 1
     \end{array}
  \right] \quad
  \left\{
     \begin{array}{lcl}
        x_1         & = & b_1(t) \\
        x_1 x_2     & = & b_2(t) \\
        x_1 x_2 x_3 & = & b_3(t)
     \end{array}
  \right.
\end{equation}
where the solutions $x_k(t)$, $k=1,2,3$, are of the form
\begin{equation}
   \exp(\alpha t) + O(t^4) = 1 + \alpha t 
   + \frac{\alpha^2}{2!} t^2
   + \frac{\alpha^3}{3!} t^3 + O(t^4),
\end{equation}
with $\alpha \in [-1, -1 + \delta] \cup [1 - \delta, 1]$, $\delta > 0$,
or $|\alpha| = 1$ for random $\alpha \in \cc$.
The $\alpha$ introduces numerical variation in the coefficients
of the solutions.

The choice of exponential series as solution series avoids that the
series in the right hand side vector of the system become large.
Considering the series expansion
\begin{equation}
   \exp(t) = \sum_{k=0}^{d-1} \frac{t^k}{k!} + O(t^d)
\end{equation}
leads in Table~\ref{tabMPneed} to a justification for multiprecision,
based on the size of the last coefficient in the truncated series.
The levels of precision proceeds in powers of two
according to the expected quadratic convergence of Newton's method.

\begin{table}[hbt]
\begin{center}
\caption{Recommended precision levels based on the order of the series,
where {\tt eps} is the smallest positive double 
that makes a difference when added to~1.0.}
\begin{tabular}{r||r|c|r}
       \multicolumn{1}{c||}{$k$} 
     & \multicolumn{1}{c|}{$1/k!$}
     & choice of precision
     & \multicolumn{1}{c}{\tt eps} \\ \hline \hline
  7  & {\tt 2.0e-004} & double precision & {\tt 2.2e-16} \\ \hline
 15  & {\tt 7.7e-013} & use double doubles & {\tt 4.9e-32} \\
 23  & {\tt 3.9e-023} & use double doubles & \\
 31  & {\tt 1.2e-034} & use quad doubles & {\tt 6.1e-64} \\
 47  & {\tt 3.9e-060} & use octo doubles & {\tt 4.6e-128} \\
 63  & {\tt 5.0e-088} & use octo doubles & \\ \hline
 95  & {\tt 9.7e-149} & need hexa doubles & {\tt 5.3e-256} \\
127  & {\tt 3.3e-214} & need hexa doubles & \\
\end{tabular}
\label{tabMPneed}
\end{center}
\end{table}

The ``{\small need hexa doubles}'' in Table~\ref{tabMPneed}
is because accelerating the least squares solving,
evaluation and differentiation at power series
with hexa double arithmetic is still a work in progress.
The computations in this paper are therefore limited to series of
order 64.  While computations of larger orders are possible,
the accuracy in octo double precision is no longer guaranteed.

While one column of monomials is sufficient for convergence
and scalability investigations, 
consider the 2-column format of monomials
\begin{equation}
   \bfc_1 \x^{E_1} + \bfc_2 \x^{E_2} = \bfb(t),
\end{equation}
for two $n$-vectors $\bfc_1$ and $\bfc_2$
and two exponent matrices $E_1$ and $E_2$.
With the introduction of new variables,
any polynomial system can be written in this 2-column format.

For the experiments with two columns of monomials, 
specific lower and upper triangular matrices of ones are used.
For example, for $n = 3$:
\begin{equation} \label{equpplow}
   E_1 =
   \left[
      \begin{array}{ccc}
         1 & 0 & 0 \\
         1 & 1 & 0 \\
         1 & 1 & 1
      \end{array}
   \right]
   \quad \mbox{and} \quad
   E_2 =
   \left[
      \begin{array}{ccc}
         1 & 1 & 1 \\
         1 & 1 & 0 \\
         1 & 0 & 0 
      \end{array}
   \right].
\end{equation}
Although the product of the degrees of the system
$\bfc_1 \x^{E_1} + \bfc_2 \x^{E_2} = \bfb(t)$ 
is now much larger than the systems defined by 
$\x^{E_1} = \bfb(t)$ and $\x^{E_2} = \bfb(t)$,
there is still only one solution,
which allows for decaying coefficients in the power series
and thus for coefficients of modest growth,
benefiting the numerical conditioning of the problem.

That $\bfc_1 \x^{E_1} + \bfc_2 \x^{E_2} = \bfb(t)$,
for lower and upper triangular matrices $E_1$
and $E_2$ as in~(\ref{equpplow}),
has only one solution can be seen by considering
the $k$-th and the $(n-k)$-th equations:
\begin{eqnarray*}
   k & : & c_{1,k\hphantom{-k}}
           \x^{{\bf e}_{1, k\hphantom{-k}}}
           + c_{2,k\hphantom{-k}}
           \x^{{\bf e}_{2, k\hphantom{-k}}} = b_k(t) \\
 n-k & : & c_{1,n-k} \x^{{\bf e}_{1, n-k}}
           + c_{2,n-k} \x^{{\bf e}_{2, n-k}} = b_{n-k}(t).
\end{eqnarray*}
As ${\bf e}_{1, k} = {\bf e}_{2, n-k}$
and ${\bf e}_{2, k} = {\bf e}_{1, n-k}$,
the two equations can be diagonalized into
\begin{eqnarray*}
   \gamma_{1, k} \x^{{\bf e}_{1, k}}
   & = & \beta_k(t) \\
   \gamma_{1, n-k} \x^{{\bf e}_{1, n-k}}
   & = & \beta_{n-k}(t),
\end{eqnarray*}
so for those two particular choices of $E_1$ and $E_2$
the system is equivalent to the one-column system
$\x^E = \bfb(t)$ which has a unique solution series.

While $\bfc_1 \x^{E_1} + \bfc_2 \x^{E_2} = \bfb(t)$
has thus the same good numerical conditioning
as $\x^E = \bfb(t)$, it serves as a good test on
the increased cost of evaluation and differentiation.
For problems with many nearby singularities,
it is recommended to work with a factor~$\delta \in (0, 1)$
to multiply the parameter~$t$ with to dampen the growth
of the coefficients in the series,
according to the numerical interpretation of Theorem~\ref{thmFabry}.

\section{Staggered Computations} \label{secStaggered}

In computing $\x(t) = \x_0 + \x_1 t + \x_2 t^2 + \cdots + \x_{d-1} t^{d-1}$,
not all $d$ coefficient vectors need to be involved.
We start $\x_0$ with half its precision correct,
otherwise Newton's method may not converge.
The first iteration consists of getting $\x_0$ correct
to the full working precision.
If Newton's method would not converge for order zero,
then there is no use of increasing the order.
The $d$ in the order $O(t^d)$ is increased gradually,
for example, the update formula for the order
\begin{equation}
   d := d + 1 + d/2
\end{equation}
is optimistically hoping for quadratic convergence.

Once $\x_k$ is correct, the corresponding $\bfb_k = 0$,
as $\bfb_k$ is obtained by evaluation, and then the update
$\Delta \x_k$ should no longer be computed because
\begin{equation}
   QR \Delta \x_k = \bfb_k = \zero
   \quad \Rightarrow \quad \Delta \x_k = \zero.
\end{equation}
This gives a criterion to stop the iterations.

\section{Accelerating Newton's Method}

While the blocked Householder QR is very suitable to GPU acceleration and
teraflop performance is achieved already at relatively modest dimensions,
it stars only at the very beginning of Newton's method as it is no longer
needed once the QR decomposition is computed.
The second part of the least squares solver, the back substitution,
is needed in every stage, as are the convolutions to compute the
right hand sides of the linear systems.

In~\cite{KH10}, the {\em Compute to Global Memory Access (CGMA) ratio}
is defined as the number of floating-point calculations performed by
a kernel for each access to the global memory.
This CGMA ratio corresponds to the more general notion of
{\em arithmetic intensity} of a computation~\cite{WWP09}.

\subsection{Arithmetic Intensity of Convolutions}

In the computation of the arithmetic intensity of convolutions,
or equivalently, the number of floating-point computations per double,
it suffices to consider one monomial.
For example, take $n=4$ and let $f = x_1 x_2 x_3 x_4$ be the monomial
we evaluate and differentiate.  Using the reverse mode of algorithmic
differentiation, the computations are organized 
in three columns as follows:
\begin{equation}
  \begin{array}{l|l|l}
     x_1 \star x_2     & x_4 \star x_3     & x_1 x_2 \star x_4 \\
     x_1 x_2 \star x_3 & x_4 x_3 \star x_2 & x_4 x_3 \star x_1 \\
     x_1 x_2 x_3 \star x_4 & &
  \end{array}
\end{equation}
where each $\star$ indicates a new multiplication.
In the three columns we count respectively $n-1 = 3$, $n-2 = 2$,  
$n-2 = 2$, for a total of $(n-1) + (n-2) + (n-2) = 3 n - 5$ multiplications,
for $n$ inputs.
If the inputs were doubles, then the arithmetic intensity would be
\begin{equation}
   \frac{3 n - 5}{n}.
\end{equation}

Each input is a power series of order~$d$.  To avoid thread divergence,
the coefficients of the second series in each product are padded with zeros.
For example, for $d = 3$:
\begin{equation}
   \begin{array}{ll}
      & (a_0 + a_1 t + a_2 t^2) (b_0 + b_1 t + b_2 t^2) \\
    = & (a_0 \star b_0 + a_1 \star b_{-1} + a_2 \star b_{-2}) \\
    + & (a_0 \star b_1 + a_1 \star b_{0\phantom{-}} + a_2 \star b_{-1}) t \\
    + & (a_0 \star b_2 + a_1 \star b_{1\phantom{-}}
         + a_2 \star b_0\phantom{-}) t^2,
   \end{array}
\end{equation}
where coefficients with negative indices are zero.
Ignoring the additions, we count $d^2$ multiplications,
so we now have $(3n-5)d^2$ multiplications for $n d$ inputs.
If the coefficients of the power series were doubles,
then the arithmetic intensity would be
\begin{equation} \label{eqaiconv}
   \frac{(3n - 5) d^2}{n d}.
\end{equation}

The coefficients of the power series are multiple doubles.
For double double, quad double, and octo double, the number of
doubles in the inputs are respectively $2 n d$, $4 n d$, and $8 n d$.
Working with complex coefficients doubles the size of the input.
Doubling the precision doubles the size of the input,
but increase the arithmetical cost significantly,
as illustrated by Table~\ref{tabcostmd} (with data from~\cite{Ver20}).

\begin{table}[hbt]
\begin{center}
\caption{Arithmetical Cost of Multiple Double Multiplications, e.g.:
multiplying two double doubles takes 5 additions of two doubles,
9 subtractions and 9 multiplications, for a total of 23
floating-point operations.}
\begin{tabular}{r||r|r|r||r}
   & \multicolumn{1}{c|}{$+$}
   & \multicolumn{1}{c|}{$-$}
   & \multicolumn{1}{c||}{$*$} & total \\ \hline
double double &   5 &   9 &   9 &   23 \\
  quad double &  99 & 164 &  73 &  336 \\
  octo double & 529 & 954 & 259 & 1742
\end{tabular}
\label{tabcostmd}
\end{center}
\end{table}

Then the number of floating-point operations per double
for the evaluation and differentiation of a product of $n$
power series of order $d$ are respectively
for double doubles, quad doubles, and octo doubles:
\begin{equation}
   \frac{23 (3n - 5) d^2}{2 n d},
   ~~ \frac{336 (3n - 5) d^2}{4 n d},
   \mbox{ and } \frac{1742 (3n - 5) d^2}{8 n d},
\end{equation}
where the corresponding multiplication factors 
$23/2$, $336/4$, and $1742/8$ 
evaluate respectively to 11.5, 84, and 217.75.

With each doubling of the precision, the arithmetic intensity
increases by a factor of 11.5 (from double to double double),
by a factor of $7.30 \approx 84/11.5$ (from double double to quad double), and
by a factor of $2.59 \approx 217.75/84$ (from quad double to octo double).
As the evaluation and differentiation is needed at every stage
of Newton's method, 
the high numbers of floating-point operations per double
are promising indicators for the success of GPU acceleration.

\subsection{Complex Vectorized Arithmetic}
The analysis in the previous section concerned real numbers.
On complex data, the computational intensity increases,
and moreover, offers more parallelism.

Writing the real and imaginary parts of $\x$ as ${\cal R}(\x)$
and ${\cal I}(\x)$, respectively, and applying the same notation to~$\y$,
the product of $\x$ with $\y$ is then executed as
\begin{equation} \label{eqcmplxvec}
\begin{array}{ccc}
         {\cal R}(\x) \star {\cal R}(\y)
   & - & {\cal I}(\x) \star {\cal I}(\y) \\
         {\cal R}(\x) \star {\cal I}(\y)
   & + & {\cal I}(\x) \star {\cal R}(\y).
\end{array} 
\end{equation}
All $\star$'s in~(\ref{eqcmplxvec}) can be executed simultaneously first,
with the $-$ and $+$ executed next.
The formulas in~(\ref{eqcmplxvec}) correspond to
what is called the 4M method in~\cite{Kou23}.

The vectorization of the complex arithmetic introduces more parallelism.
Alternatively, the computation can be arranged such that one block of
threads simultaneously collaborate on the real and imaginary parts
of one product.  This alternative is beneficial for vector shorter
than the typical warp size of 32 threads.

\subsection{Accelerated Least Squares}

With the QR decomposition of~$A$,
solving $A \x = \bfb$ in the least squares sense
is reduced to $R \x = Q^H \bfb$,
to the multiplication of $Q^H$ with $\bfb$,
followed by a back substitution.

Assuming the leading coefficient vector $\x_0$ of the series
has an accuracy of at least half the working precision,
the QR decomposition happens only once in the first stage
of Newton's method, at a cost of $O(n^3)$.
While every stage involves the solution of $R \x = Q^H \bfb$,
the cost of computing $Q^H \bfb$ and the back substitution
is both $O(n^2)$.  Even if the acceleration of the QR decomposition
works better than the acceleration of $R \x = Q^H \bfb$,
the factor $n$ in the cost overhead of QR over $R \x = Q^H \bfb$
makes that the proportion of the QR decomposition will still
dominate all $Q^H \bfb$ computations and all back substitutions,
as $n$ equals 1024 and the number of stages is capped to~24.

\subsection{Accelerated Updates and Residuals}

The updates of the right hand side vectors
involve many different matrices.
For example, the updates to $\bfb_3$ happen as
$\bfb_3 := \bfb_3 - A_3 \x_0$, $\bfb_3 := \bfb_3 - A_2 \x_1$,
$\bfb_3 := \bfb_3 - A_1 \x_2$,
each time with different matrices $A_3$, $A_2$, and $A_1$
which cannot remain all in the main memory of the device.
Even as the cost of these computations is $O(n^2)$
we may expect the updates to occupy a significant portion
of the total execution times.

The same arguments apply to the computation of the residuals,
when measuring the accuracy of the computed updates $\Delta \x$.
However, one could significantly reduce the cost by selecting
only one or a few random equations instead of computing the
residuals for all equations.

For comparison with the arithmetic intensities of the convolutions,
consider the matrix-vector product, for an $n$-by-$n$ matrix.
Executed on doubles, $n^2$ multiplications are performed on
$(n+1) n$ doubles.  The ratio $n^2/(n^2 + n)$ improves on
power series of order~$d$, as convolutions with padding take
$d^2$ multiplications, while the size of the data is multiplied by~$d$.
Restricting to multiplications, the arithmetic intensity then is
\begin{equation}
   \frac{n^2 d^2}{(n+1)n d} \approx d.
\end{equation}
Compared to the arithmetic intensity of convolutions,
the leading terms of the numerator of~(\ref{eqaiconv})
divided by the denominator evaluates to $3 d$.
So the arithmetic intensity of convolutions is three times more 
than that of matrix-vector products.

\section{Staging Multiple Doubles}

This section details 
the implementation of the multiple double arithmetic.

\subsection{Inlining of Arithmetical Kernels}
The obvious storage of a multiple double
is as an array of many doubles, stored at consecutive locations.
But this storage is unfavorable for coalesced memory accesses 
when working with vectors of multiple doubles.
Therefore, all arithmetical operations on multiple doubles
are defined by functions which take as input arguments not
arrays of doubles, but sequences of doubles.
For example, the function to add two double doubles has prototype
\begin{verbatim}
     void ddf_add
      ( double a_hi, double a_lo,
        double b_hi, double b_lo,
        double *c_hi, double *c_lo );
\end{verbatim}
where high doubles of the numbers $a$, $b$, and $c$ are 
respectively given as 
{\tt a\_hi}, {\tt b\_hi}, {\tt c\_hi}
and their low doubles as
{\tt a\_lo}, {\tt b\_lo}, {\tt c\_lo},
for the addition $c := a + b$.
The GPU counterpart to the {\tt ddf\_add} above has the same arguments
but starts with
\begin{verbatim}
  __device__ __forceinline__  void ddg_add
\end{verbatim}
forcing the inlining of the kernels.
All kernels (in the {\tt .cu} files) are compiled with {\tt nvcc -O3}.

Every supported precision level has dedicated functions,
once defined as regular C functions (for correctness comparisons)
and once defined as kernel functions to be executed on the device.

\subsection{Shared Memory and Registers}
A vector of double doubles is stored in main memory on the device
as two vectors, one vector with the high doubles, and another vector
with the low doubles.  All threads in one block load the high and
low doubles in vectors in memory shared by all threads in the block.
As adjacent threads read into adjacent memory locations,
this staging the multiple doubles benefits memory coalescing.

The inlined arithmetical kernels are applied to numbers in registers.
Before executing the arithmetical operations, the data is copied
from shared memory into local variables, variables that are local
to each thread in a block.  The additional instruction cycles
allow to hide the latencies of the data transfers.

\section{Computational Results}

The computational experiments attempt to answer the
following three questions.
Classifying the types of kernels into three categories:
convolutions, least squares, updates and residuals,
which type of kernel occupies the largest portion of
the overall execution time?
For which order of the series do we reach teraflop performance?
What happens to the wall clock time when the precision is doubled?

The starting point used for performance comparisons, or the baseline,
is the observed performance compared against the relative performance
of the GPUs.  For example, 
if one GPU is theoretically twice as fast as another,
then we may expect that computations on the faster GPU end in half
the time of the other.

\subsection{Graphics Processing Units}

The code was developed for the ``Volta'' V100 NVIDIA GPU,
and tested on the ``Pascal'' P100 and RTX 2080, 4080 NVIDIA GPUs.
The ``Ampere'' A100 was added most recently.
Table~\ref{tabgpus} lists the main characteristics of the GPUs.

\begin{table}[hbt]
\begin{center}
\caption{Specifications of the GPUs.  
All have 64 cores per multiprocessor (MP),
except for the RTX 4080, which has 128 cores per~MP.}
\begin{tabular}{r||r|r|r|r||r|r|r|r}
  \multicolumn{1}{c||}{GPU}  
&  CUDA  & \#MP  & 
  \#cores & GHz &
  \multicolumn{2}{c|}{host CPU} & RAM & GHz \\ \hline
      Pascal P100 & 6.0~~ &  56~~ & 
3584~ & 1.33 &
 \multicolumn{2}{c|}{Intel E5-2699} & 256 GB &  2.20 \\
      Volta V100 & 7.0~~ &  80~~ &  
5120~ & 1.91 &
 \multicolumn{2}{c|}{Intel W2123} & 32 GB & 3.60 \\
     Ampere A100 & 8.0~~ & 108~~ &
6912~ & 1.41 &
 \multicolumn{2}{c|}{Intel 5318Y} & 256 GB & 2.10 \\ \hline
  RTX 2080 & 7.5~~ &  46~~ &  
2944~ & 1.10 &
 \multicolumn{2}{c|}{Intel i9-9880H} & 32 GB & 2.30 \\
  RTX 4080 & 8.9~~ &  58~~ & 7424~ & 2.10 &
 \multicolumn{2}{c|}{Intel i9-13900HX} & 16 GB & 2.20
\end{tabular}
\label{tabgpus}
\end{center}
\end{table}

The double precision peak performance of the P100 is 4.7 TFLOPS.
At 7.8 TFLOPS, the V100 is 1.66 times faster than the P100.
The theoretical peak performance of the A100 is 8.7~TFLOPS.
To evaluate the algorithms, compare the ratios of the wall clock times
on the P100 over V100 with the factor 1.66.
For every kernel, the number of arithmetical operations is accumulated.
The total number of double precision operations is computed
using the cost overhead multipliers.

The double precision peak performance of the RTX 2080 Max Q-design
is 201.5 GFLOPS, which corresponds to 761.5 GFLOPS for
the more recent RTX~4080.
Both the RTX~2080 and~4080 are housed in a laptop computer.


The operating system of the host of the P100 and V100 is
CentOS~7.  For the A100, the host runs Rocky Linux.
The code on the host is compiled with 
{\tt g++}~4.8.5~(Red Hat), with optimization flag {\tt '-O3'},
using version~8.0 of the CUDA Software Development Kit.
The RTX~2080 and RTX~4080 are housed in laptops running
Windows~11.  The code on the host of RTX~2080 is compiled
with the 2019 Community Visual Studio compiler,
using CUDA~10.2.  On the host of RTX~4080, the compiler
is the 2022 Community Visual Studio version, using
CUDA~12.1.  On Windows, the optimization flag for the
compiler is~{\tt '-02'}.




\subsection{Proportions of Kernel Times}

Table~\ref{tab6kernels} lists six different types of kernels.

\begin{table}[hbt]
\begin{center}
\caption{Six different kernels.}
\begin{tabular}{c|c}
name & legend \\ \hline
convolutions & evaluation and differentiation \\
qr & Householder QR \\
qhb &  $Q^H \bfb$ computations \\
bs & back substitutions for $R \x = Q^H \bfb$ \\
updates & updates $\bfb := \bfb - A \x$ \\
residuals & residual computations $\| \bfb - A \x \|_1$
\end{tabular}
\label{tab6kernels}
\end{center}
\end{table}

\begin{figure}[hbt]
\centerline{\includegraphics[width=9.0cm,
trim = 1cm 2.5cm 0cm 0cm, clip]{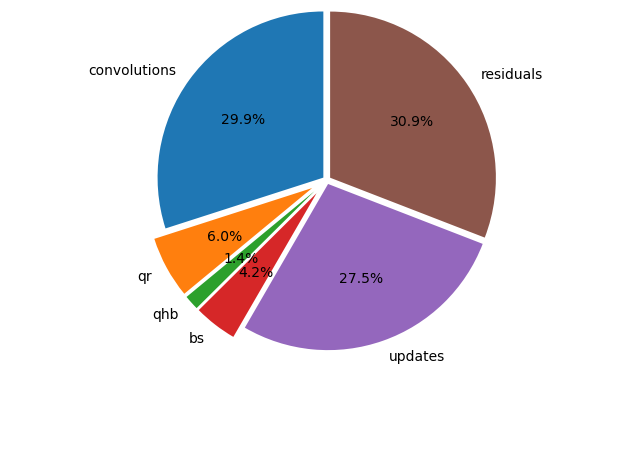}}

\caption{Percentage of each type of accelerated computation 
for a one column monomial system in octo double precision, on V100,
with legend in Table~\ref{tab6kernels}.}
\label{figpiechart}
\end{figure}

Visualizing the data in Table~\ref{tabpietimes},
Figure~\ref{figpiechart} shows the percentages of the kernels 
on one column of monomials defined by a triangular 
exponent matrix of dimension 1024 to compute series
of order 64 in octo double precision, done on the V100.

\begin{table}[hbt]
\begin{center}
\caption{Times in seconds for each kernel on V100,
 on one column and two columns of monomials in octo double precision.}
\begin{tabular}{r|r|r}
                 kernel & one column & two columns \\ \hline
            convolution &  121.386~~ &  244.535~~ \\
         Householder QR &   24.451~~ &   24.123~~ \\
  kernel for $Q^H \bfb$ &    5.849~~ &    6.139~~ \\
      back substitution &   17.053~~ &   17.884~~ \\
                updates &  111.474~~ &  124.080~~ \\
              residuals &  125.122~~ &  137.963~~ \\ \hline
     total kernel times &  405.334~~ &  554.723~~ \\ 
        wall clock time & 1129.794~~ & 1808.480~~
\end{tabular}
\label{tabpietimes}
\end{center}
\end{table}

\begin{figure}[hbt]
\centerline{\includegraphics[width=9.0cm,
trim = 1.3cm 1.2cm 1.3cm 1.7cm, clip]{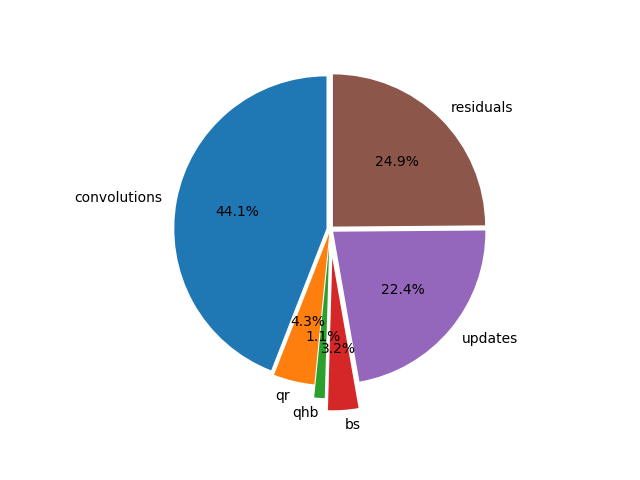}}

\caption{Percentage of each type of accelerated computation 
for a two column monomial system in octo double precision, on V100,
with legend in Table~\ref{tab6kernels}.}
\label{figpiechart2}
\end{figure}

The largest portion of the time goes to the residual computations,
for all equations in the system.
The residuals are important to measure the convergence
and must be computed in multiple double precision.
One optimization could be to select at random one or a couple
of equations and compute the residuals for those selected equations
instead of for all equations.
Figure~\ref{figpiechart2} shows that for a 2-column monomial system,
the time spent on convolutions dominates.

\subsection{Performance of Convolutions}

Figure~\ref{figpiechart} shows that the convolutions occupy a substantial
part.  For what orders of the series do we
observe teraflop performance?  Consider Table~\ref{tabperfcnv}.

\begin{table}[hbt]
\begin{center}
\caption{Performance in gigaflops of convolutions on P100, V100, and A100
  to evaluate and differentiate one column of 1024 monomials
  in octo double precision.}
\begin{tabular}{c|r|r|r}
  order & \multicolumn{1}{c|}{P100} & 
          \multicolumn{1}{c|}{V100} &
          \multicolumn{1}{c}{A100} \\ \hline
 1 &    8.041  &   28.997 & 39.284 \\
 2 &   16.191  &   59.820 & 95.770 \\
 3 &   23.748  &   90.003 & 147.440 \\
 5 &   29.277  &  149.894 & 252.502 \\
 8 &   62.747  &  240.760 & 411.946 \\
12 &   94.035  &  360.816 & 623.255 \\
18 &  140.918  &  540.572 & 938.634  \\
27 &  211.261  &  810.645 & 1412.459 \\
41 &  351.994  & 1045.032 & 1601.007 \\
62 &  535.136  & 1569.347 & 2472.086 \\
64 &  554.654  & 1658.382 & 2568.016
\end{tabular}
\label{tabperfcnv}
\end{center}
\end{table}

\begin{figure}[hbt]
\begin{center}
\begin{picture}(350,210)
\put(-25,110){\includegraphics[width=14cm]{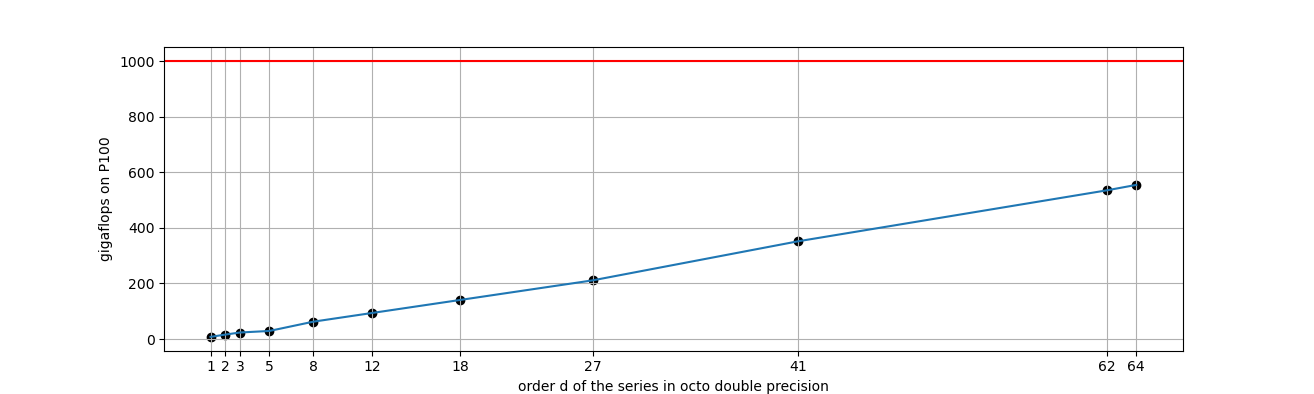}}
\put(-25, 0){\includegraphics[width=14cm]{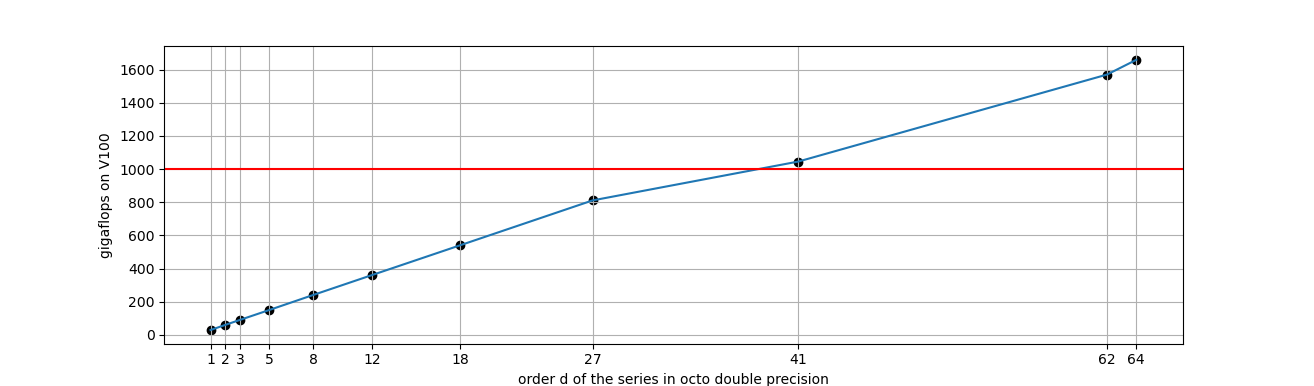}}
\end{picture}
\caption{Performance in gigaflops on the P100 (top) and on the V100 (bottom)
to evaluate and differentiation at series in octo double precision
versus the order of the series.}
\label{figconvolutions}
\end{center}
\end{figure}

On one column of monomials, triangular exponent matrix of ones,
$n=1024$, performance of the evaluation and differentiation,
in octo double precision, for increasing orders of the series,
visualizing the data in Table~\ref{tabperfcnv},
Figure~\ref{figconvolutions} shows that teraflop performance is
observed after order 40 on the V100.
Teraflop is observed after order 40 on the V100,
but on the P100 only half a teraflop is reached at order~64.
On the A100, teraflop performance is observed at order~27.

In the implementation, 
in the convolution of two power series
each thread is responsible for one coefficient of the result.
Threads are launched in blocks of the size that match
the number of coefficients. 
An implementation better suited for series of lower order
would employ a finer granularity and have several threads
collaborate to compute one coefficient of a convolution of two series.

As the V100 is theoretically 1.66 faster than the P100,
multiplying 554.654 (on the last row of Table~\ref{tabperfcnv})
by 1.66 yields 920.726, which is less than the observed 1,658.382 
on the V100.  Thus, for higher levels of precision (octo double),
newer GPUs scale better than older ones.
The observed 1.7 teraflop performance on the V100 is only at 22\%
of the theoretical peak performance.
The performance improves on series of higher order.

\subsection{Doubling the Precisions}

To investigate
how much of the cost overhead can be compensated by the acceleration
consider the wall clock times and the elapsed times spend by all kernels
when the precision is doubled.

\begin{table}[hbt]
\begin{center}
\caption{Wall clock times in seconds
for double (D), double double (2D), quad double (4D),
and octo double (8D) for one column of 1024 monomials 
to compute series of order 64 on P100 and V100.}
\begin{tabular}{r|r|r|r|r}
   & \multicolumn{1}{c|}{D} 
   & \multicolumn{1}{c|}{2D} 
   & \multicolumn{1}{c|}{4D} 
   & \multicolumn{1}{c}{8D} \\ \hline
P100 kernel times &  10.4 &  44.5 & 204.4 &  1263.3 \\
     wall clock   & 418.5 & 695.3 & 969.9 &  3073.6 \\ \hline
V100 kernel times &   6.2 &  22.6 & 146.4 &   405.3 \\ 
     wall clock   & 277.9 & 475.3 & 834.7 &  1129.8
\end{tabular}
\label{tabwall}
\end{center}
\end{table}

\begin{figure}[hbt]
\begin{center}
\begin{picture}(340,360)(30,0)
\put(0,180){\includegraphics[width=14cm]{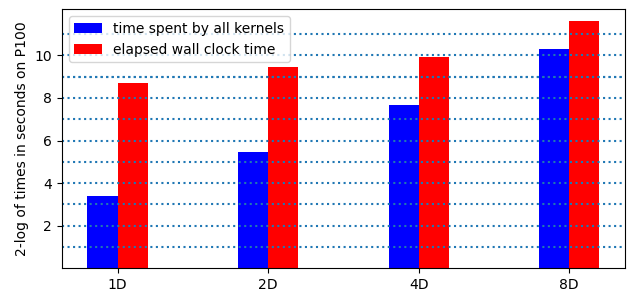}}
\put(0, 0){\includegraphics[width=14cm]{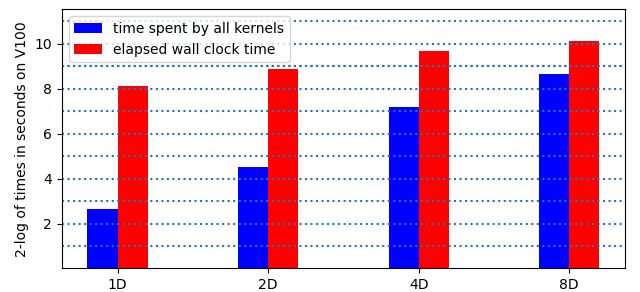}}
\end{picture}
\caption{In doubling the precision, 
the wall clock times on P100 (top) and V100 (bottom) less than double
as the proportion of the elapsed times spent by all kernels increases,
on one column of~1,024 monomials.}
\label{figwall}
\end{center}
\end{figure}

Visualizing Table~\ref{tabwall},
Figure~\ref{figwall} shows the 2-logarithms of the times
of 24 steps with Newton's method on one column of monomials
defined by a triangular exponent matrix of ones of dimension 1024,
on the V100.
Doubling the precision less than doubles the wall clock time
and increases the time spent by all kernels.
The reason for this is that problems become more computationally intensive:
the amount of data doubles but the amount of arithmetical computations 
on average is multiplied by factors larger than two,
as explained in section~5.1, see also Table~\ref{tabcostmd}.

\subsection{Experiments on the RTX 2080 and 4080}

\begin{table}[hbt]
\begin{center}
\caption{Wall clock times in seconds
for double (D), double double (2D), quad double (4D),
and octo double (8D) for one column of 512 monomials 
to compute series of order 64 on RTX~2080.}
\begin{tabular}{r|r|r|r|r}
   & \multicolumn{1}{c|}{D} 
   & \multicolumn{1}{c|}{2D} 
   & \multicolumn{1}{c|}{4D} 
   & \multicolumn{1}{c}{8D} \\ \hline
kernel times &  1.4 &  16.6 & 122.4 & 380.1 \\
  wall clock & 35.7 &  80.0 & 225.0 & 474.8
\end{tabular}
\label{tabwallrtx2080}
\end{center}
\end{table}

On the problem of applying the block Householder QR on real matrices
of octo double numbers, the RTX~4080 reaches a performance
of 413.2 GFLOPS on dimension 1,024 (8 blocks of size 128),
and 719.5 GFLOPS (very close to the theoretical peak performance)
is reached on dimension 2,048 ($8~\times~256$).

\begin{figure}[hbt]
\centerline{\includegraphics[width=14cm]{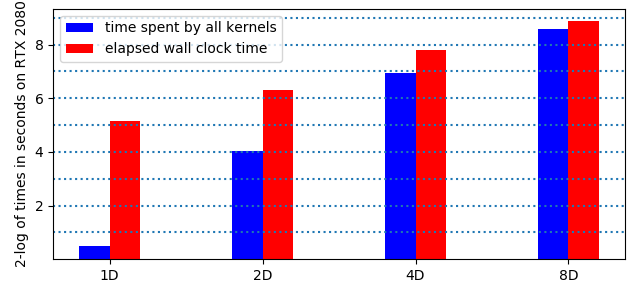}}
\caption{Doubling the precision on the RTX 2080,
on one column of 512 monomials.}
\label{figwallrtx2080}
\end{figure}

The last experiments concern the RTX 2080,
on 16 steps with Newton's method on one column of monomials
defined by a triangular exponent matrix of ones of dimension 512.
The results are summarized in Figure~\ref{figwallrtx2080},
visualizing the data in Table~\ref{tabwallrtx2080}.

Table~\ref{tabrtx4080vs2080} compares the performance
of the newer RTX 4080 versus the older RTX 2080,
illustrating the scalability, on one column of 512 monomials,
on series of order~64.
Although the RTX 4080 is at least twice as powerful as the RTX 2080,
the order~64 is unfavorable for the RTX 4080 which has 128 cores
per multiprocessor.
Relative to the P100 and V100, Figure~\ref{figwallrtx2080}
looks better than the corresponding Figure~\ref{figwall}
because the RTX 2080 is already better occupied,
already at a column of~512 monomials.

\begin{table}[hbt]
\begin{center}
\caption{Kernel times in seconds and gigaflops performance
for the Householder QR, multiplication of $Q$ with $\bfb$ (Qb),
back substitution (BS), evaluation \& differentation (ED).}
\begin{tabular}{c|rr|rr}
 & \multicolumn{2}{c|}{RTX 2080}
 & \multicolumn{2}{c}{RTX 4080} \\
stage &   time & gigaflops &   time & gigaflops \\ \hline
 QR   &   28.5 &  100.5~~ &   18.3 &  156.5~~ \\
 Qb   &   18.6 &   26.3~~ &   14.9 &   33.0~~ \\
 BS   &   18.1 &   16.3~~ &   16.4 &   17.2~~ \\
 ED   &   27.3 &  110.0~~ &   15.8 &  190.1~~
\end{tabular}
\label{tabrtx4080vs2080}
\end{center}
\end{table}

\section{Conclusions}

Although several optimizations in the code will improve the
performance, this first implementation offers a promising first
step towards a scalable nonlinear solver based on results from
numerical analytic continuation.

Using decaying coefficients of power series expansions,
octo double precision suffices for series of order~64.
Teraflop performance of the evaluation and differentiation
is already attained at order 40 on the V100.
The convolutions to evaluate and differentiate at power series
remain a significant portion of all computational work.
For two columns of monomials which can encode general
polynomial systems, the computational effort to evaluate and
differentiate dominates. 

Doubling precisions less than doubles the wall clock times
because the computations are then compute bound and thus
well suited for acceleration by graphics processing units.
Extending the acceleration beyond octo double precision on 
GPUs is a future direction.

\medskip
\noindent {\bf Acknowledgements.}  Preliminary versions of this study were
presented at the 2023 Joint Mathematical Meetings, 
in the AMS Special session on Polynomial Systems, Homotopy Continuation
and Applications (organized by Tim Duff and Maggie Regan); and
in the minisymposium on Co-Design for Heterogeneous System Architectures
(organized by Navamita Ray and Edgar Leon) at the SIAM CSE 2023 meeting.
The author thanks the organizers of those meetings.
%
%
%
\bibliographystyle{plain}

\end{document}